\documentclass{article}

\usepackage{latexsym}
\usepackage{epsfig}
\usepackage{amssymb}
\usepackage{amsthm}
\usepackage{amsmath}
\usepackage{epstopdf}

\def\1{{\bf 1}}

\newtheorem{thm}{Theorem}[section]
\newtheorem{cor}[thm]{Corollary}
\newtheorem{lem}[thm]{Lemma}

\newtheorem{rem}[thm]{Remark}
\newtheorem{defi}[thm]{Definition}

\newenvironment{prf}{\noindent{\bf Proof:} }{\hfill$\Box$\mbox{}}

\title{The Kite Graph is Determined by Its Adjacency Spectrum}

\author{
Hatice Topcu\thanks{\small Dept. of Mathematics,
Nev\c{s}ehir Hac{\i} Bekta\c{s} Veli University, Turkey}\\
\small{\tt{haticekamittopcu@gmail.com}}
\\[5pt]
Sezer Sorgun\footnotemark[1]\\
\small{\tt{srgnrzs@gmail.com}}
}
\date{}

\begin{document}

\maketitle

\begin{abstract}

\noindent
 The Kite graph $Kite_{p}^{q}$ is obtained by appending the complete graph $K_{p}$ to a pendant vertex of the path $P_{q}$. In this paper, the kite graph is proved to be determined by the spectrum of its adjacency matrix.
\end{abstract}
AMS subject classification 05C50.

\section{Introduction}
In this paper, all graphs are simple and undirected. Let $G=(V(G),E(G))$ be a graph with vertex set $V(G)=\{v_{1}, v_{2}, \ldots ,v_{n}\}$ and edge set $E(G)$. The \textit{order of G} is the number of the vertices in $G$. The \textit{degree} of a vertex $v$ is denoted by $d(v)$. The least degree in $G$ is denoted by $\delta(G)$ and the largest by $\Delta(G)$. For $i\in \{1, \ldots, n\}$, $N(v_{i})$ denotes the set of the vertices of $G$ which are adjacent to $v_{i}$. An edge contains a vertex of degree $1$ is called a \textit{pendant edge}.  A $clique$ of $G$ is a complete subgraph of $G$ and \textit{size of a clique} is the number of the vertices in the clique. Number of the vertices in the largest clique of $G$ is called the $clique$ $number$ of $G$ and denoted by $w(G)$. A clique with size $3$ in $G$ is called a \textit{triangle } in $G$. A \textit{petal} is added to a graph when we add a pendant edge and then duplicate this edge to form a pendant $2$-cycle. A \textit{blossom} $B_{k}$ consists of $k$ petals ($k\geq0$) attached a single vertex; thus $B_{0}$ is just the trivial graph. A graph with blossoms (possibly empty) at each vertex is called a $B-graph$.  The $line$ $graph$ $L(G)$ of a graph $G$ is the graph whose vertices are the edges of $G$, with two vertices in $L(G)$ adjacent whenever the corresponding edges in $G$ have exactly one vertex in common. Let $H$ be a graph with the vertex set $\{h_{1}, \ldots, h_{n}\}$ and let $a_{1}, \ldots, a_{n}$ be non-negative integers. The \textit{generalized line graph} $GL(H; a_{1}, \ldots, a_{n})$ is the graph $L(\widehat{H})$, where $\widehat{H}$ is the $B-graph$ $H(a_{1}, \ldots, a_{n})$ obtained from $H$ by adding $a_{i}$ petals at vertex $h_{i}$ $(i=1, \ldots, n)$.

Let $A(G)$  be the $(0,1)-adjacency$ $matrix$ $of$ $G$ and $P_{A(G)}(\lambda)=det(\lambda I-A(G))$ is the $characteristic$ $polynomial$ $of$ $G$ where $I$ is the $identity$ $matrix$. Since the matrix $A(G)$ real and symmetric, its eigenvalues,i.e. all roots of $P_{A(G)}(\lambda)$, are all real numbers and called the $adjacency$ $eigenvalues$ $of$ $G$. We use $\lambda_{1}(G)\geq\lambda_{2}(G)\geq\ldots\geq\lambda_{n}(G)$ to denote the adjacency eigenvalues of $G$. These eigenvalues compose the $adjacency$ $spectrum$ of $G$. The largest adjacency eigenvalue of $G$ is known as its $adjacency$ $spectral$ $radius$ and denoted by $\rho(G)$. Let $D(G)$ be the diagonal matrix of vertex degrees of $G$. Then the matrices $\xi(G)=D(G)-A(G)$ and $Q(G)=D(G)+A(G)$ are called  the \textit{Laplacian matrix} and the \textit{signless Laplacian matrix} of \textit{G}, respectively. The \textit{Laplacian spectrum} and the \textit{signless Laplacian spectrum} consists of the eigenvalues of $\xi(G)$ and \textit{Q(G)}, respectively. We use $\mu_{1}(G)\geq\mu_{2}(G)\geq\ldots\geq\mu_{n}(G)=0$ to denote the signless Laplacian eigenvalues of $G$.  Two graphs $G$ and $H$ are said to be $cospectral$ (with respect to the adjacency, Laplacian or signless Laplacian matrix,etc.) if they share the same (adjacency, Laplacian or signless Laplacian, etc.) spectrum. A graph $G$ is said to be \textit{determined by its spectrum} if there isn't any other non-isomorphic graph that is cospectral with $G$.

 If a tree has exactly one vertex of degree greater than 2, then it is called a \textit{starlike tree}. A starlike tree with maximum degree $\Delta$ is denoted by $T(l_{1}, l_{2}, \ldots, l_{\Delta})$ such that $T(l_{1}, l_{2}, \ldots, l_{\Delta})- v = P_{l_{1}} \cup P_{l_{2}} \cup \ldots \cup P_{l_{\Delta}}$ where $v$ is the vertex of degree $\Delta$ in the starlike tree. The \textit{Kite graph}, denoted by $Kite_{p}^{q}$, is obtained by appending a complete graph $K_{p}$ to a pendant vertex of a path $P_{q}$. A $Kite_{p}^{q}$ graph is actually the line graph of a starlike tree $T(l_{1}, l_{2}, \ldots, l_{p})$ where $l_{1}=l_{2}=\ldots=l_{p-1}=1$ and $l_{p}=q+1$. The line graph of any starlike tree is called a \textit{sunlike graph} in some paper\cite{10}. So, a kite graph is also a sunlike graph with some conditions. The \textit{Double kite graph}, denoted by $DK(p,q)$, is obtained by appending  complete graph $K_{p}$ to both pendant vertices of a path $P_{q}$ \cite{27}.

The Kite graph already appeared in the literature several times \cite{1,10,20,24,25,26,27}. In \cite{1}, authors have shown that among connected graphs with order $n$ and clique number $w$, the minimum value of the adjacency spectral radius is attained for a Kite graph. Also, they've given a small interval for the adjacency spectral radius of the Kite graph with any fixed clique number. Nath and Paul\cite{24} have shown that among all graphs, the Kite graph is the unique graph that maximizes the distance spectral radius. Among all connected graphs with a fixed clique number, the first four smallest values of the adjacency spectral radius and the first three smallest values of the Laplacian spectral radius are obtained in \cite{25} and \cite{26}, respectively. By the help of all of these papers, it can be said that the Kite graph plays an important role on the extremum values of the spectral radius of some graph matrices.

Also, Stani$\grave{c}$ \cite{27} has proved that, in connected graphs any Double kite graph is determined by its adjacency spectrum. By appending a cycle to a pendant vertex of a path, we obtain the well-known lollipop graph. It is already shown that the lollipop graphs are determined by their spectrum \cite{3,13}. Similarly, by appending a complete graph to a pendant vertex of a path, we obtain the Kite graph. In \cite{20}, the adjacency characteristic polynomial of the Kite graph is calculated and it is shown that the $Kite_{p}^{2}$ graph is determined by its adjacency spectrum. Also it is conjectured that the $Kite_{p}^{q}$ graph is determined by its adjacency spectrum for all $p$ and $q$. By the motivation of the all of these results, in this paper we prove that the  $Kite_{p}^{q}$ graph is determined by its adjacency spectrum for all $p$ and $q$.

Additionally, the problem : "Which non-regular graphs with least eigenvalue at least -2 are determined by their adjacency spectrum?" is proposed firstly by van Dam and Haemers \cite{5} and partially answered in some recent results \cite{14,19}. Zhou and Bu have proved that the line graph of a starlike tree with maximum degree $\Delta$ is determined by its adjacency spectrum when $\Delta\geq12$\cite{19}. They have left an open problem for $\Delta<12$:"When is the line graph of a starlike tree determined by its adjacency spectrum?". Since the smallest eigenvalue of the Kite graph is greater than -2 and the $Kite_{p}^{q}$ graph is actually the line graph of a starlike tree, our result in this paper partially answers the open problems of Dam-Haemers\cite{5} and Zhou-Bu\cite{5}.

  For $p\leq2$, $p=3$ and $q=0$, clearly we obtain the path graph, the lollipop graph and  the complete graph, respectively. It is well-known that these graphs are already determined by their adjacency spectrum \cite{5}. Hence, we will consider the case $p\geq4$.


\section{Preliminaries}
\label{}
\begin{lem}\cite{5,19}
Adjacency eigenvalues of the Path graph with $n$ vertices are 2cos$\frac{\pi j}{n+1} (j=1, \ldots, n)$.

\end{lem}

\begin{lem} \cite{8}\textbf{(Interlacing Lemma)}

If G is a graph on n vertices with eigenvalues $\lambda_{1}(G)\geq\ldots\geq\lambda_{n}(G)$ and H is an induced subgraph on m vertices with eigenvalues $\lambda_{1}(H)\geq\ldots\geq\lambda_{m}(H)$, then for $i=1, \ldots , m$ \[\lambda_{i}(G)\geq\lambda_{i}(H)\geq\lambda_{n-m+i}(G) \]

\end{lem}

\begin{lem}\cite{19,21}
Let G be a graph with n vertices and m edges, and let L(G) be the line graph of G. Then

\[P_{A(L(G))}(x) = (x+2)^{m-n}P_{Q(G)}(x+2) \]

\end{lem}

For a graph $G$, the $subdivision$ $graph$ of $G$, denoted by $S(G)$, is the graph obtained from $G$ by inserting a new vertex in each edge of $G$.

\begin{lem}\cite{19,21}
Let G be a graph with n vertices and m edges, and let S(G) be the subdivision graph of G. Then

\[P_{A(S(G))}(x) = x^{m-n}P_{Q(G)}(x^{2}) \]

\end{lem}

\begin{thm}\cite{18}
Each starlike tree with maximum degree 4 is determined by its signless Laplacian spectrum.
\end{thm}

\begin{lem}\cite{21}
The starlike tree $T(l_{1}, \ldots, l_{\Delta}) (\Delta\geq5)$ is determined by its signless Laplacian spectrum.

\end{lem}

A connected graph whose largest adjacency eigenvalue is equal to 2 is called a \textit{Smith graph}. All Smith graphs are shown in Figure 1.

\begin{lem}\cite{8}
A connected graph with the largest adjacency eigenvalue less than 2 are precisely the induced proper subgraphs of the Smith graphs shown in Figure 1.

\end{lem}

\begin{center}\begin{figure}[h]
 \includegraphics[height=7cm,width=14cm]{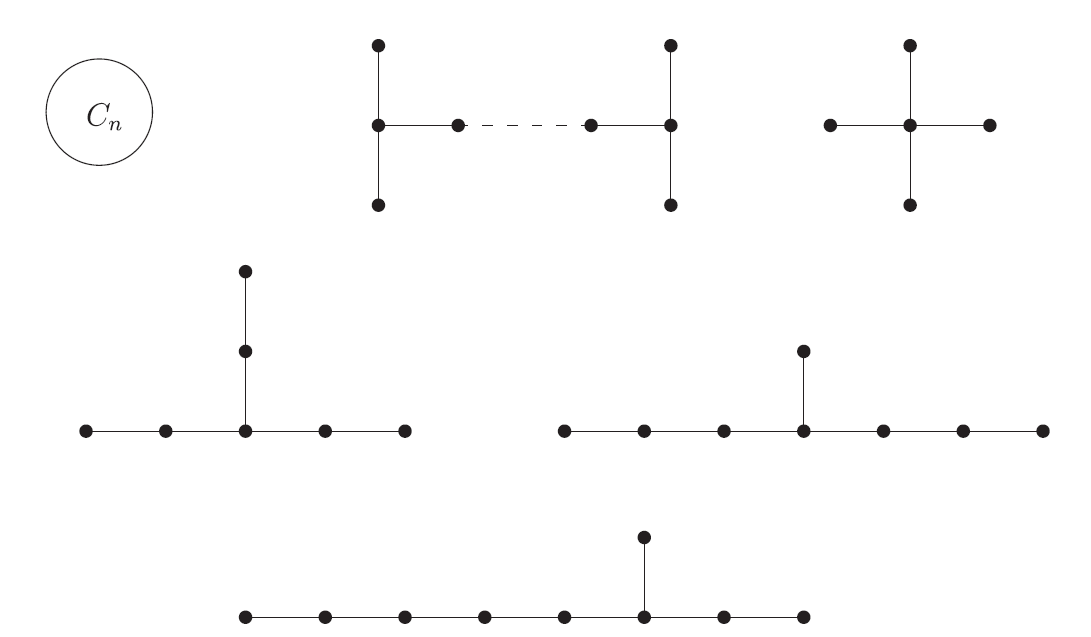}\caption{ The Smith Graphs}\label{fig:Smith Graph}
\end{figure}\end{center}

An \textit{exceptional graph} is a connected graph with the least adjacency eigenvalue at least -2 which is not a generalized line graph. All connected graphs with the least adjacency eigenvalue greater than -2 are characterized by Doob and Cvetkovi\'{c} \cite{22}.

\begin{lem}\cite{14,22}
Let G be a connected graph with n vertices. Then $\lambda_{n}(G)>-2$ if and only if one of the following holds:

(i) G = L(H) where H is a tree or an odd unicyclic graph.

(ii) G = GL(H;1,0,$\ldots$,0) where H is a tree.

(iii) G is one of the 573 exceptional graphs.

\end{lem}

All exceptional graphs and their adjacency spectral radiuses are listed in [9,Appendix Table A2].

\begin{defi}\cite{14,19}
Let G be a graph with n vertices. The quantity \[\prod_{i=1}^{n} (\lambda_{i}(G)+2)\] is called the discriminant of G and denoted by $d_{G}$.
\end{defi}

If $G$ has more than one component, then the discriminant of $G$  is equal to the product of discriminants of all components of $G$. If two graphs have the same adjacency spectrum, then they have same discriminant.

\begin{lem}\cite{14}
Let G be a connected graph with least adjacency eigenvalue greater than -2. The following statements hold :

(i) If G is an exceptional graph with eight vertices, then $d_{G}$=1.

(ii) If G is an exceptional graph with seven vertices, then $d_{G}$=2.

(iii) If G is an exceptional graph with six vertices, then $d_{G}$=3.

(iv) If G is the line graph of an odd unicyclic graph, then $d_{G}$=4.

(v) If $G = GL(\widehat{H})$ where $\widehat{H}$ is the B-graph obtained from a tree by attaching a single petal, then $d_{G}$=4.

(vi) If G is the line graph of a tree and G has n vertices, then $d_{G}=n$+1.

\end{lem}

\begin{lem} \cite{19}
If two graphs G and H have the same signless Laplacian spectrum, then their line graphs have the same adjacency spectrum. The converse is true if G and H have the same number of vertices and edges.

\end{lem}

\begin{thm} \cite{19}
Let $G=T(l_{1}, \ldots, l_{\Delta})$ with $\Delta\geq12$. Then $L(G)$ is determined by its adjacency spectrum.

\end{thm}

\begin{lem} \cite{8}
Let x$_{1}$ be a pendant vertex of a graph G and x$_{2}$ be the vertex which is adjacent to x$_{1}$. Let G$_{1}$ be the induced subgraph obtained from G by deleting the vertex x$_{1}$. If x$_{1}$ and x$_{2}$ are deleted, the induced subgraph G$_{2}$ is obtained. Then,

\[P_{A(G)}(x)=x P_{A(G_{1})}(x)-P_{A(G_{2})}(x)\]

\end{lem}

\begin{lem} \cite{5} For $n\times n$ matrices $A$ and $B$, followings are equivalent :

\textbf{(i)} $A$ and $B$ are cospectral

\textbf{(ii) }$A$ and $B$ have the same characteristic polynomial

\textbf{(iii) }$tr(A^{i})=tr(B^{i})$ for $i=1,2,...,n$

\end{lem}

\begin{lem} \cite{5} For the adjacency matrix of a graph $G$, the following parameters can be deduced from the spectrum;

\textbf{(i)} the number of vertices

\textbf{(ii)} the number of edges

\textbf{(iii)} the number of closed walks of any fixed length.

 \end{lem}

\begin{lem} \cite{3} Let $H_{n,p}$ denotes the graph that is obtained by appending a cycle $C_{p}$ to a pendant vertex of a path $P_{n-p}$. Then, the graph $H_{n,p}$ with $p$ odd is determined by its adjacency spectrum.
\end{lem}

\begin{lem} \cite{9} A graph is a generalized line graph if and only if it has no induced subgraph isomorphic to one of the graphs $G^{1}, \ldots, G^{31}$ from Figure 3.
\begin{defi}If a graph $G$ can not contain $H$ as an induced subgraph, then $H$ is called a \textit{\textbf{forbidden subgraph}} for $G$.\end{defi}

\begin{center}\begin{figure}[h]
 \includegraphics[height=139 mm,width=14cm]{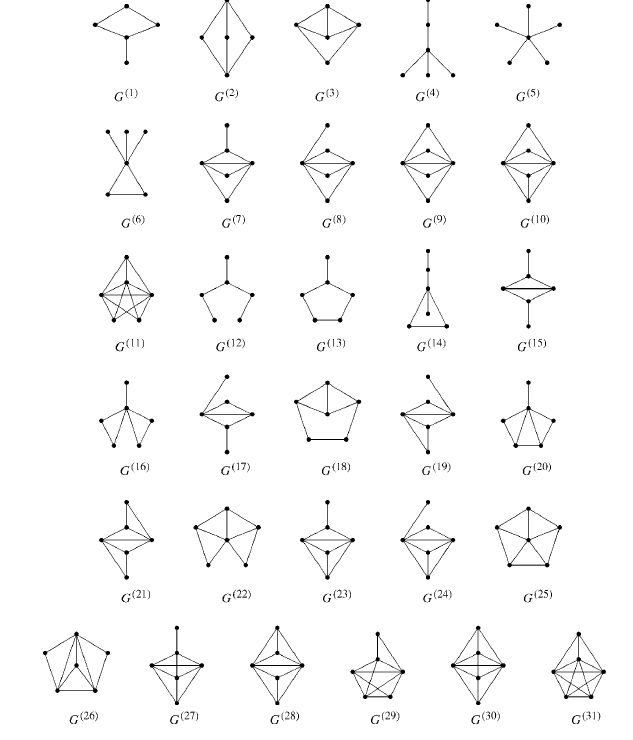}
\caption{The Forbidden subgraphs for a generalized line graph}\label{The Forbidden subgraphs for a generalized line graph}\end{figure}
 \end{center}

\end{lem}

\section{Main Results}
\label{}

Firstly, we will give some lemmas and corollaries which will be used to prove the main result of this paper.  The following lemma could be obtained from the proof of Theorem 2.2.

\begin{lem} \cite{19}
Let G be a starlike tree and $\lambda_{2}(L(G))$ denotes the second largest adjacency eigenvalue of the line graph of G. Then $\lambda_{2}(L(G))<$2.
\end{lem}

\begin{prf}
Let $G=T(l_{1}, l_{2}, \ldots, l_{\Delta})$ and $m=l_{1}+ l_{2}+ \ldots + l_{\Delta}$. By Lemma 2.4,  adjacency eigenvalues of the subdivision graph $S(G)$ are $\pm\sqrt{\mu_{1}(G)},$  $\pm\sqrt{\mu_{2}(G)},$  $\ldots,$  $\pm\sqrt{\mu_{m}(G)}, 0$. Since $S(G)$ is also a starlike tree, we can assume that $v$ is the vertex of degree $\Delta$ in $S(G)$. Hence,  $S(G)-v = P_{2l_{1}} \cup P_{2l_{2}} \cup \ldots \cup P_{2l_{\Delta}}$. By Lemma 2.1 and Lemma 2.2, we have $\sqrt{\mu_{2}(G)}<2$ and  so $\mu_{2}(G)<4$. Then, by Lemma 2.3, $\lambda_{2}(L(G))<2$.
\end{prf}

\begin{lem} \cite{19}
Let $G=T(l_{1}, l_{2}, \ldots, l_{\Delta})$, $\Delta\neq3$ and H is a connected graph. If H is cospectral with the L(G) w.r.t. adjacency spectrum, then H must be isomorphic to the L(G).
\end{lem}

\begin{prf}
Since $G$ is a starlike tree, we know that $\Delta\geq2$. If $\Delta=2$, then $G$ is a path and so the $L(G)$ is a path. It is well-known that the path graphs are determined by their adjacency spectrum \cite{5}. By Lemma 2.12, $L(G)$ and $H$ have the same number of vertices and edges. If $\Delta\geq4$, then the order of the $L(G)$ is greater than 4. Also by Lemma 2.8, we have $disc(L(G))=n+1$. From these facts, we get $disc(H)=disc(L(G))>5$. So, by Lemma 2.7 and Lemma 2.8, $H$ is the line graph of a tree, $\widetilde{T}$.  Since $L(G)=L(T(l_{1}, l_{2}, \ldots, l_{\Delta}))$ and $H=L(\widetilde{T})$, we obtain that $G$ and $\widetilde{T}$ have the same number of vertices and edges. Lemma 2.9 implies that $G$ and $\widetilde{T}$ have the same signless Laplacian spectrum. From Theorem 2.1 and Lemma 2.5, we can say that $H$ is isomorphic to the $L(G)$. \end{prf}

\begin{rem}

For the case $\Delta=3$, we refer \cite{3,12}.

\end{rem}

Following result could be obtained directly from Lemma 3.2.

\begin{cor}
Among connected graphs, the line graph of each starlike tree with maximum degree different from 3 is determined by its adjacency spectrum.
\end{cor}

\begin{cor}
Among connected graphs, the $Kite_{p}^{q}$ graph is determined by its adjacency spectrum for all $p$ and $q$.
\end{cor}

\begin{prf}
This result could be obtained automatically from Lemma 2.13 and Corollary 3.1.

\end{prf}

\begin{lem}
If any given graph $G$ is adjacency cospectral with the $Kite_{p}^{q}$ graph, then $G$ must be connected.

\end{lem}

\begin{prf}

If $p\leq2$ or $p=3$, we obtain the path graph or the lollipop graph. Hence, we continue with $p\geq4$. A kite graph actually the line graph of a starlike tree such that $Kite_{p}^{q}=L(T(l_{1}, \ldots , l_{p}))$ where $l_{1}= \ldots = l_{p-1}=1$ and $l_{p}=q+1$. By Theorem 2.2, the $Kite_{p}^{q}$ graph is determined by its adjacency spectrum for $p\geq12$. Hence, we continue with $4\leq p\leq 11$. From Lemma 3.1, we have $\lambda_{2}(Kite_{p}^{q})<2$. By Lemma 2.7, we get $\lambda_{n}(Kite_{p}^{q})>-2$ such that $n=p+q$.

Assume that a graph $G$ is cospectral with the $Kite_{p}^{q}$ with respect to the adjacency spectrum and $G$ is disconnected. Let $G = H \cup H_{r_{1}} \cup \ldots \cup H_{r_{k}}$ (k$\geq$1) where $H, H_{r_{i}}$ $ (i=1, \ldots, k)$ are the components of $G$ and $\lambda_{1}(G) = \lambda_{1}(H)$. From here and the interlacing lemma, we get $\lambda_{2}(H)\leq \lambda_{2}(G)<2$ and $\lambda_{1}(H_{r_{i}})\leq \lambda_{2}(G)<2$ $(i = 1, \ldots, k)$. By Lemma 2.6, $H_{r_{i}}$ $(i=1, \ldots, k)$ are the induced proper subgraphs of the Smith graphs. So that, the components of $G$ apart from $H$ are trees and we get $t(G)=t(H)$. Since $G$ is cospectral with the $Kite_{p}^{q}$, by Lemma 2.12,  we have \[t(G)=t(Kite_{p}^{q})=(\begin{array}{c}
                                                                              p \\
                                                                              3
                                                                            \end{array}) = t(H)\]

By calculation, we have  $\lambda_{1}(Kite_{p}^{q})<p$  for the largest eigenvalue of the $Kite_{p}^{q}$. As it can be seen easily, the clique number of $H$, \[w(H)\leq p\] because $\lambda_{1}(H) = \lambda_{1}(Kite_{p}^{q})<p$.

\vspace{0.5 cm}
\textsl{\textbf{Claim 1:}} $w(H)<p$.
\vspace{0.5 cm}

Let $w(H) = p$. If $n_{H}$ denotes the order of $H$, then $n_{H}<n$. Since $t(H) = (\begin{array}{c}
                                                                              p \\
                                                                              3
                                                                            \end{array}) $, there could be just one clique with the size $p$ and there is no more triangle out of this clique. We already know that $\lambda_{2}(H)<2$ and $\lambda_{n_{H}}(H)>-2$. Hence, from the interlacing lemma,  $H$ can not contain any cycle (except of the clique) and the star graph with four branches as induced subgraphs. Also, the following two graphs in Figure 3 can not be induced subgraphs of $H$ because $\lambda_{1}(K_{p}^{2})>\lambda_{1}(H)$ and $\lambda_{1}(Urchin_{p}^{2})>\lambda_{1}(H)$.

\vspace{-0.7 cm}
 \begin{center}\begin{figure}[h]
 \includegraphics[height=5.5cm,width=13cm]{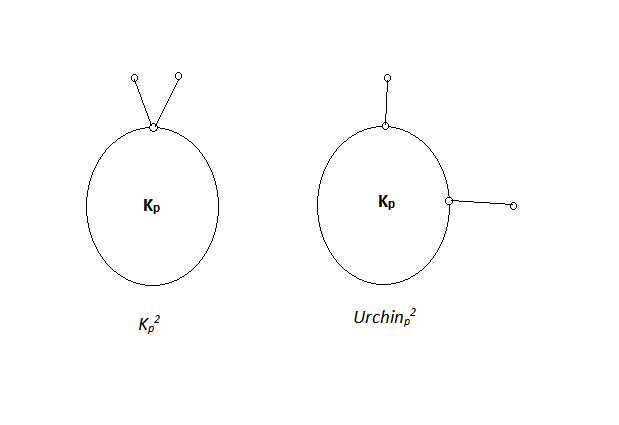}
\caption{The graphs K$_{p}^{2}$ and Urchin$_{p}^{2}$}\label{K$_{p}^{2}$ and Urchin$_{p}^{2}$}\end{figure}
 \end{center}

\vspace{-0.7 cm}

Thus, $H$ could be just in the form that is the graph obtained by attaching a tree to a vertex of the clique $K_{p}$. Let us denote this tree by $T$ . Since the $S_{4}$ is a forbidden subgraph for $H$, the degree of any vertex in $T$ could be equal to at most 3. In Figure 4, $T'$ denotes a tree which has exactly two vertices with degree 3 and other vertices that are between these two vertices have degree 2.  Also, the least and the largest eigenvalues of $T'$ are equal to -2 and 2, respectively. Hence, $H$ can not contain $T'$ as an induced subgraph. Actually, $T'$ is a Smith graph. So that, $T$ may have at most one vertex with degree 3.
\vspace{-0.5 cm}
\begin{center}\begin{figure}[h]
 \includegraphics[height=35mm,width=11cm]{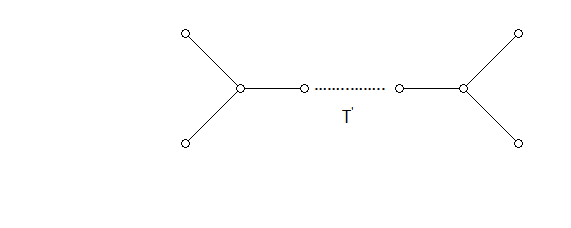}
\caption{The tree $T'$}\label{The tree $T'$}\end{figure}
 \end{center}
\vspace{-0.5 cm}
\newpage
Let us label this vertex with $v$ in $T$ such that $d(v)=3$. If $v$ is adjacent to an endpoint of $T$, then $0\in spec(H)$. But, we've already known that $0\not\in spec(Kite_{p}^{q})$. Hence $v$ can not be adjacent to an endpoint of $T$. Now, we consider the other positions of $v$ in $T$ by the following figure.

\begin{center}\begin{figure}[h]
 \includegraphics[height=9cm,width=14cm]{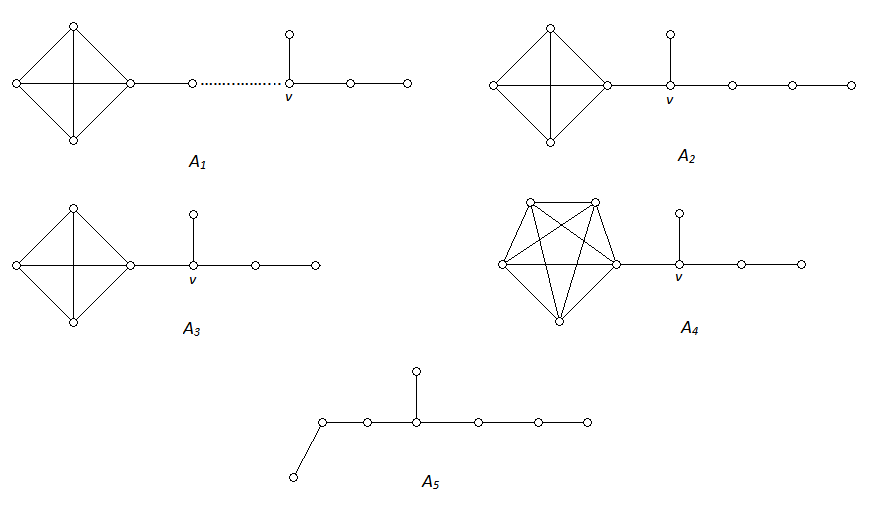}
\caption{The graphs A$_{1}$, A$_{2}$, A$_{3}$, A$_{4}$, A$_{5}$ }\label{The graphs A$_{1}$, A$_{2}$, A$_{3}$, A$_{4}$, A$_{5}$}\end{figure}
 \end{center}

If the orders of $A_{1}, A_{2}, A_{3}, A_{4}, A_{5}$ denoted by $n_{1}, n_{2}, n_{3}, n_{4}, n_{5}$, respectively, then $\lambda_{n_{1}}(A_{1})\leq-2$, $\lambda_{n_{2}}(A_{2})\leq-2$ and $\lambda_{n_{4}}(A_{4})=\lambda_{n_{5}}(A_{5})=-2$. So, $H$ does not contain $A_{1}, A_{2}, A_{4}, A_{5}$ as induced subgraphs. Hence, the situation that is the vertex $v$ is contained by $T$ is just possible in the case that $H$ is isomorphic to $A_{3}$. But in this case, for the largest eigenvalues, $\lambda_{1}(A_{3})>\lambda_{1}(Kite_{p}^{q})$. Thus, $T$ can not contain any vertex with degree 3. This means that, all of the vertices in $T$ have degree at most 2 and so $T$ is a path. If $q'$ denotes the order of $T$, then $H$ is isomorphic to the $Kite_{p}^{q'}$. When $q'<q$, we've already known that $spec(Kite_{p}^{q'})\not\subseteq spec(Kite_{p}^{q})$. This fact contradicts with $spec(H)\subseteq spec(Kite_{p}^{q})$. Accordingly, we have $w(H)\neq p$. Hence, we showed that our first claim is true.

\vspace{0.3 cm}

So, we get $3\leq w(H)\leq10$ since $4\leq p \leq 11$ and $w(H)<p$.

\vspace{0.3 cm}

Since $\lambda_{1}(H)=\lambda_{1}(Kite_{p}^{q})$ and $4 \leq p \leq 11$, by the calculation and  the comparison, we can say that $H$ is not an exceptional graph [9, Appendix, TableA2]. From this fact, and Lemma 2.7, $H$ is a generalized line graph of a tree with a single petal attached or a line graph of a tree (or an odd unicyclic graph). Therefore, showing the truth of the following claim will be enough to complete the proof.

\vspace{0.5 cm}
\textsl{\textbf{Claim 2:}} $H$ can not be a generalized line graph of any graph.
\vspace{0.5 cm}

In this step, we first give some notations.

\vspace{0.3 cm}

For $i\in V$ and $V\subseteq V(H)$, $T_{i}$ denotes the set of the adjacent vertices to $i$ such that these vertices are not adjacent to any other vertex in $V$. Hence, if $x\in T_{i}$, we get $x\sim i$ and $x\nsim y$ for all $y\in V-\{i\}$. For $i,j\in V$, $T_{ij}$ denotes the set of the adjacent vertices to both $i$ and $j$ at the same time such that any vertex in this set is not adjacent to any other vertex in $V$. So, if $x\in T_{ij}$, then we get $x\sim i$, $x\sim j$ and $x\nsim y$ for all $y\in V-\{i,j\}$. More generally, for all $x\in T_{i_{1}\ldots i_{k}}$ and $y\in V-\{i_{1}\ldots i_{k}\}$, we obtain $x\sim i_{1}$, $\ldots$, $x\sim i_{k}$ and $x\nsim y$ such that $|V|=n'$, $k\in \{1, \ldots, n'\}$ and $i_{k}\in V$.

Graphs in Figure 6 are forbidden subgraphs for $H$ because of their second or least adjacency eigenvalues. Also, for $i\in\{1,2,3,4,5,6\}$, $n(F)$ and $n(F_{i})$ denotes the orders of $F$ and $F_{i}$, respectively.

\begin{center}\begin{figure}[h]
 \includegraphics[height=130mm,width=130mm]{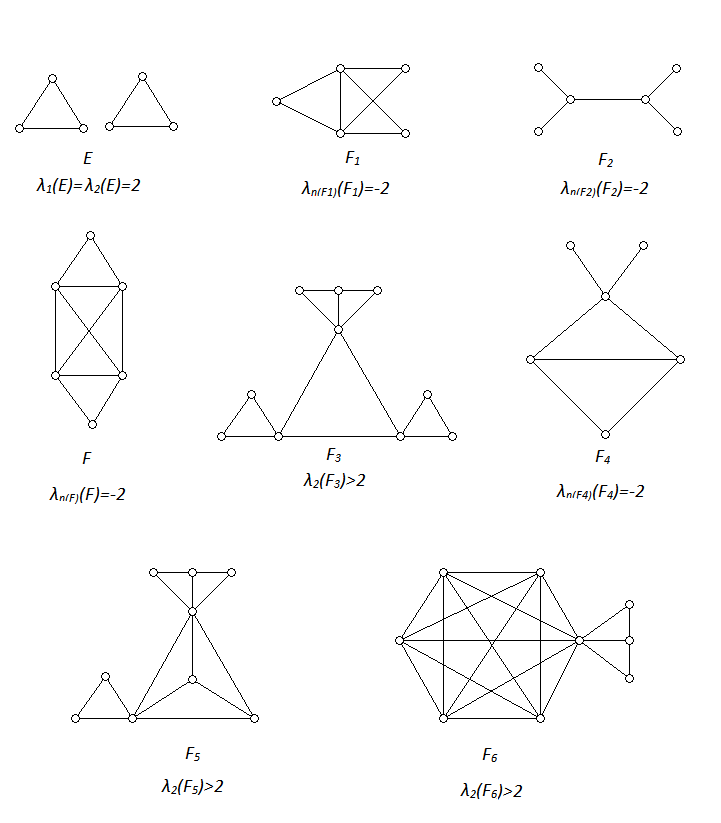}
\caption{The graphs E, F, F$_{1}$, F$_{2}$, F$_{3}$, F$_{4}$, F$_{5}$ and F$_{6}$ }\label{}\end{figure}
 \end{center}

Let $H$ be a generalized line graph. By Lemma 2.14, the graphs that are given in Figure 2 are forbidden subgraphs in this case. The graphs from Figure 7 are also used along the rest of the proof.

\begin{center}\begin{figure}[h]
 \includegraphics[height=140mm,width=133mm]{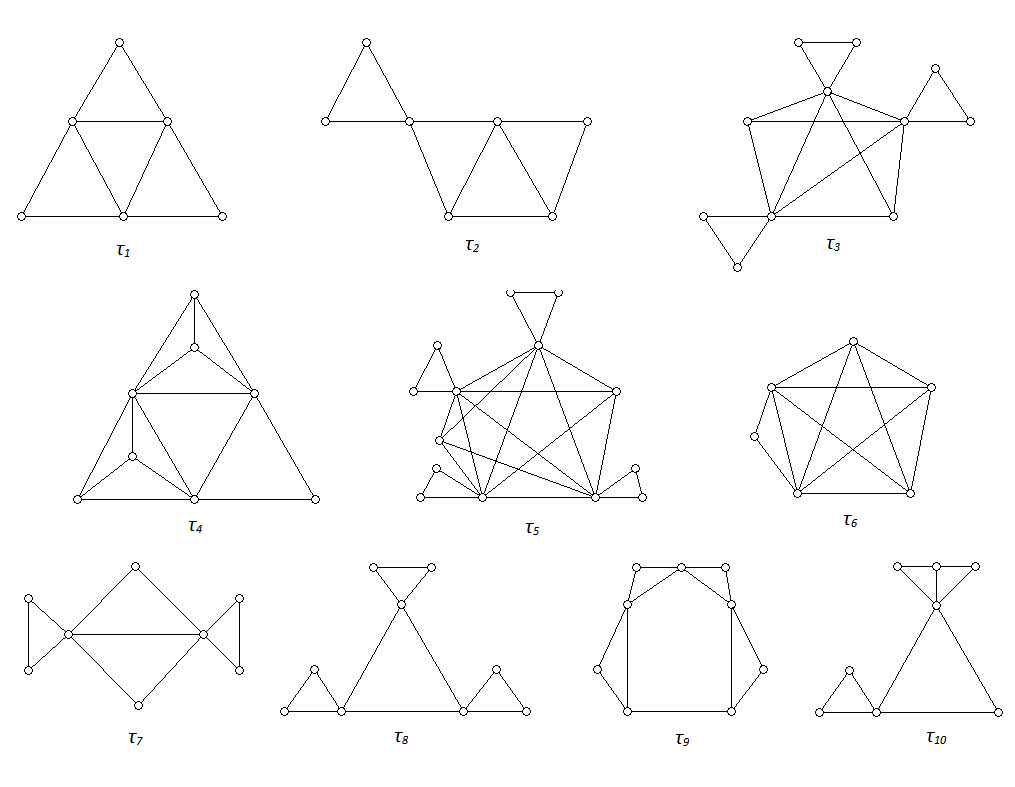}
\caption{The graphs $\tau_{1}, \tau_{2}, \tau_{3}, \tau_{4}, \tau_{5}, \tau_{6}$, $\tau_{7}$, $\tau_{8}$, $\tau_{9}$, $\tau_{10}$ }\label{ The graphs} \end{figure}
 \end{center}

We will consider the situation, that $H$ is a generalized line graph, in four cases.

\textbf{Case 1 :} Let $w(H)=3$. If we label the vertices of a triangle of $H$ with $V'=\{1,2,3\}$, then it is clear that $T_{123}=\emptyset$. $\forall i,j\in V', x,y\in T_{ij}$, if $x\sim y$ or $x\nsim y$ then $K_{4}$ or $F_{1}$ will be an induced subgraph,respectively. Thus we get $\forall i,j\in V', |T_{ij}|\leq1$. Now, we consider the sets $T_{i}$, for all $i\in V'$. There could be at most two edges between any three vertices of $T_{i}$ to prevent the occurrence of $K_{4}$. Moreover, we can say that the number of the edges between any three vertex in $T_{i}$ can not be equal to 0 or 1 because $G^{16}$ and $S_{4}$(or $G^{6}$) are forbidden subgraphs. Therefore, by using the fact that $G^{20}$ is a forbidden subgraph, we get $\forall i \in V', |T_{i}|\leq3$. Simultaneously, since  $F_{2}$ and $F_{3}$ are forbidden subgraphs, we get, if $ \forall i\in V', T_{i}\neq\emptyset$ then $\sum_{i\in V'}|T_{i}|\leq6$ and if $ \exists i\in V', T_{i}\neq\emptyset$ then $\sum_{i\in V'}|T_{i}|\leq5$. Hence, we continue with the following subcases.

\textbf{Subcase 1 :} Assume that for all $i,j\in V'$ at least two of the $T_{ij}$'s are non-empty. Without loss of generality, let $T_{12}=\{x\}$ and $T_{13}=\{y\}$. If $x\sim y$ then $C_{4}$ will be an induced subgraph of $H$. Hence, $x\nsim y$. By using the forbidden subgraphs $C_{4}, K_{4}, G^{19}, G^{20}, G^{21},$ $G^{22}, G^{24}, G^{25}$, we can conclude that $|T_{x}|\leq2, |T_{y}|\leq2$ and $|T_{x}|+|T_{y}|\leq3$. Also, all of the other sets $T_{1}, T_{2}, T_{3}, T_{xy},\ldots, T_{xy123}$ are empty. Since $E$ is a forbidden subgraph, any triangle in $H$ must contain at least one vertex (or edge) of $\overset{\triangle}{123}$ or be adjacent to the $\overset{\triangle}{123}$. Thereby, under these conditions, the maximum value of the triangle number of $H$ could be equal to at most 5. As we said before, $t(Kite_{p}^{q})\geq10$ for $p>4$. If $p=4$,then $t(H)=4$. So that, $\tau_{1}$ or $\tau_{2}$ will be an induced subgraph of $H$. This is a contradiction since $\rho(Kite_{4}^{q})<\rho(\tau_{1})<\rho(\tau_{2})$.

\textbf{Subcase 2 :}Assume that just one of the $T_{ij}$'s is non-empty. Without loss of generality, let $T_{12}=\{x\}$ and $T_{13}=T_{23}=\emptyset$. If $T_{x1}$ or $T_{x2}$ is non-empty, then this case is actually equivalent to the case that at least two of the $T_{ij}$'s are non-empty. So, we continue for $T_{x1}=T_{x2}=\emptyset$. Since $C_{4}$ and $K_{4}$ are forbidden subgraphs, we obtain that $T_{x123}=T_{123}=T_{x23}=T_{x13}=T_{x12}=T_{x3}$. Also, $\forall i\in \{x,1,2,3\}, |T_{i}|\leq2$ because $S_{4}$ and $F_{4}$ are forbidden. For all $i,j\in \{x,1,2,3\}$ and $i\neq j$, if $\exists z_{1}, z_{2}\in V(H), z_{1}\sim i$ and $z_{2}\sim j$, then $z_{1}\nsim z_{2}$. Otherwise, $F_{1}$ or $C_{4}$ will be an induced subgraph of $H$. Since $G^{15}$ and $G^{17}$ are forbidden subgraphs, we get $[|T_{1}|+|T_{2}|+|T_{3}|+|T_{x}|]\leq4$. From these facts and $E$ is a forbidden subgraph, maximum value of the $t(H)$ could be equal to at most 4. Since $t(H)=t(Kite_{p}^{q})\geq4$, we get $t(H)=4$. But, in this case, $\tau_{7}$ will be an induced subgraph of $H$. This result contradicts with $\rho(\tau_{7})>\rho(Kite_{4}^{q})$.

\textbf{Subcase 3 :}Assume that all of the $T_{ij}$'s are empty. It means that, any two of the triangles in $H$ can not share an edge. Since $E$ is a forbidden subgraph, any of the triangles in $H$ must contain or be adjacent to the at least one vertex of $V'=\{1,2,3\}$. By using the forbidden subgraphs again, we conclude the following results; $\forall i\in V, |T_{i}|\leq3$; $\forall i\in V$, if $T_{i}\neq\emptyset$ then $\sum_{i\in V}|T_{i}|\leq6$; $\exists i\in V$, if $T_{i}=\emptyset$ then $\sum_{i\in V}|T_{i}|\leq5$; $max\{t(H)\}=5$.

We have already known that, $t(H)\geq10$ for $p>4$. Also, if $p=4$ then $t(H)=4$. This implies that; $\tau_{8}, \tau_{9}$ or $\tau_{10}$ will be an induced subgraph of $H$. But $\tau_{8}$ and $\tau_{9}$ can not be induced subgraphs of $H$. If we add just one pendant edge to the $\tau_{10}$ and denote this graph with $\tau_{11}$, then this will contradict with the spectral radiuses of $\tau_{11}$ and $Kite_{4}^{q}$ such that $\rho(\tau_{11})>\rho(Kite_{4}^{q})$. Hence, we obtain that $w(H)\neq3$.

\textbf{Case 2 :} Assume that $w(H)=4$. If we label the vertices of a clique with size 4 in $H$ with $V''=\{1,2,3,4\}$, then we get $T_{1234}=\emptyset$.

Without loss of generality, let $T_{124}=\{x\}$. Since $G^{27}$ is forbidden,  we get $T_{x}=\emptyset$. By using the fact that $G^{29}$, $G^{30}$, $C_{4}$ are forbidden subgraphs and $w(H)=4$, we can say that $N(x)=\{1,2,4\}$. Also, we get $\forall i, j \in V''$, $T_{ij}=\emptyset$  since $F$ and $G^{29}$ are forbidden subgraphs. If we consider the sets $T_{i}$ for all $i \in V''$, we have $T_{3}=\emptyset$ and $|T_{1}|\leq2$, $|T_{2}|\leq2$, $|T_{4}|\leq2$ because $G^{27}$, $F_{4}$ and $E$ are forbidden subgraphs. Hence, the maximum value of the number of the triangles in $H$, $t(H)$ could be equal to at most 10 such that $|T_{1}|=|T_{2}|=|T_{4}|=2$. Since $w(H)=4$, we get $p\geq5$ and $t(Kite_{p}^{q})=t(H)\geq\left(
                                            \begin{array}{c}
                                              5 \\
                                              3 \\
                                            \end{array}
                                          \right)=10
$. So, we have $t(H)=10$ and $\tau_{3}$ will be an induced subgraph of $H$. But this is a contradiction with $\rho(\tau_{3})>\rho(Kite_{5}^{q})$. Thus we say that $\forall i,j,k \in V''$, $T_{ijk}=\emptyset$.

Without loss of generality, let $T_{12}=\{x\}$. Since $G^{23}$ and $G^{24}$ are forbidden subgraphs, we get $T_{1}=T_{2}=T_{x}=\emptyset$. Also, $E, F_{1}$ and $F_{4}$ are forbidden subgraphs, so we get $|T_{3}|\leq2$, $|T_{4}|\leq2$, $|T_{x1}|\leq2$, $|T_{x2}|\leq2$. If equality holds, then two vertices in the same set must be adjacent. At the same time, for all $i,j,k\in V''$, $T_{x3}$, all of the sets $T_{x4}$, $T_{xij}$, $T_{xijk}$, $T_{x1234}$ are empty because $C_{4}$, $F_{1}$, $G^{28}$ are forbidden subgraphs, $w(H)=4$ and $T_{ij}$'s are already empty apart from $T_{12}$. If $|T_{3}|=|T_{4}|=2$, then we get $|T_{1x}|\leq1$ and $|T_{2x}|\leq1$. So, $t(H)$ will be equal to 9. But we already know that $t(H)\geq10$ for $w(H)=4$ and $w(H)<p$. So, we have $|T_{3}|\leq1$ and $|T_{4}|\leq1$. From here, the maximum value of the $t(H)$ could be equal to at most 13 where $|T_{1x}|=|T_{2x}|=2$. Hence, we get $t(H)=10$ and $p=5$ such that $|T_{1x}|=2$, $|T_{2x}|=1$ or $|T_{1x}|=1$, $|T_{2x}|=2$. In this case $\tau_{4}$ must be an induced subgraph of $H$ which contradicts with $\rho(\tau_{4})>\rho(Kite_{5}^{q})$. Thus, we can say that $\forall i,j\in V'', T_{ij}=\emptyset$. For $i\in V''$, suppose that at least one of the sets $T_{i}$ is non-empty. Since $F_{5}$ and $E$ are forbidden subgraphs, if $\exists i\in V'', |T_{i}|=3$, then $\forall j\in V''-\{i\}$, $|T_{j}|\leq1$. Hence, $t(H)$ will be the maximum when $\forall i\in V'', |T_{i}|=2$. So we get $t(H)=\left(
        \begin{array}{c}
          4 \\
          3 \\
        \end{array}
      \right)+4
$. But this result contradicts with $t(H)\geq10$. Thus, we have $w(H)\neq4$.

\textbf{Case 3 :} Assume that $w(H)=5$. Let us label the vertices of a clique with size 5 in $H$ with $V'''=\{1,2,3,4,5\}$. Clearly $T_{12345}=\emptyset$. Since $G^{26}, G^{29}, G^{30}, G^{31}$ and $F$ are forbidden subgraphs, we get the following results for all $i,j,k,l \in V'''$ ; $\sum_{i,j\in V'''} |T_{ij}|\leq1$, $T_{ijk}=\emptyset$, $\sum_{i,j,k,l\in V'''} |T_{ijkl}|\leq1$, $\sum_{i,j,k,l\in V'''} |T_{ijkl}|+|T_{ij}|\leq1$ and $|T_{i}|\leq3$. Also if $\exists i\in V''', |T_{i}|=3$, then there must be two edges between the three vertices in the same set $T_{i}$. Without loss of generality, let $T_{1234}=\{x\}$. Then, for all $i,j,k \in V'''$, $T_{ijk}=T_{ij}=\emptyset$. Since $G^{27}$ is forbidden, we get $T_{x}=T_{5}=\emptyset$ and $H$ does not contain any vertex which is adjacent to the both $x$ and (some or all of the vertices of) $V'''$. Also, $\forall i\in V''', T_{ix}=\emptyset$ because $G^{29}$ and $C_{4}$ are forbidden subgraphs. Moreover, $\forall i,j,k,l \in V'''$, $T_{xij}=T_{xijk}=T_{xijkl}=\emptyset$ because $G^{31}, C_{4}$ are forbidden subgraphs and $w(H)=5$. At the same time, we get $\forall i \in \{1,2,3,4\}$, $|T_{i}|\leq2$ and if equality holds, then the two vertices in the same set $T_{i}$ must be adjacent. Since $w(H)=5$ and $w(H)<p$, $t(Kite_{p}^{q})=t(H)=\geq\left(
                                                         \begin{array}{c}
                                                           6 \\
                                                           3 \\
                                                         \end{array}
                                                       \right)=20
$. In view of these facts, the maximum value of the $t(H)$ could be equal to 20 and so $p=6$. But in this case, $\tau_{5}$ will be an induced subgraph of $H$ and this contradicts with the fact that $\rho(\tau_{5})>\rho(Kite_{6}^{q})$. So we get, $\forall i,j,k,l \in V'''$, $T_{ijkl}=\emptyset$. Without loss of generality, let $T_{23}=\{x\}$. Since $G^{23}, G^{24}, G^{31}, F_{4}, C_{4}, G^{28}$ are forbidden subgraphs and $w(H)=5$, we get the following results. $\forall i,j,k,l \in V'''$, $T_{x}=T_{2}=T_{3}=\emptyset$, $T_{x1}=T_{x4}=T_{x5}=\emptyset$, $T_{xij}=T_{xijk}=T_{xijkl}=\emptyset$, $T_{x12345}=\emptyset$. $|T_{x2}|\leq1$, $|T_{x3}|\leq1$, $|T_{1}|\leq2$, $|T_{4}|\leq2$, $|T_{5}|\leq2$. Also for the sets $T_{1}, T_{4}, T_{5}$, if equality holds then the two vertices in the same set must be adjacent. Additionally, $E$ is  a forbidden subgraph. Thus, the maximum value of the $t(H)$ could be equal to 13. This means that $\tau_{6}$ is an induced subgraph of $H$. But this result contradicts with the fact that $t(H)\leq20$. So we get $\forall i,j\in V'''$, $T_{ij}=\emptyset$.

Let us consider the sets $\forall i \in V'''$, $T_{i}$. If $\exists i\in V''', |T_{i}|=3$ then $\forall j\in V'''-\{i\}, |T_{j}|\leq1$ because $F_{5}$ and $E$ are forbidden subgraphs. Hence, the maximum value of the $t(H)$ could be obtain when $|T_{i}|=2$ for all $i\in V'''$ and this value is equal to $\left(
                     \begin{array}{c}
                       5 \\
                       3 \\
                     \end{array}
                   \right)+5=15
$. But, we have already known that $t(H)\geq20$. So this is a contradiction. Thus, we have $w(H)\neq5$.

\textbf{Case 4 :} Assume that $w(H)\geq6$. $K_{w}$ denotes a clique of $H$ which has maximum size and $W$ denotes the set of the vertices of this clique. Let us label the four vertices of $K_{w}$ with $V''=\{1,2,3,4\}\subset W$. Since $G^{31}$ is forbidden, we get $T_{1234}\subset K_{w}$ and $T_{ijk}=\emptyset$. $\forall i,j,k \in V''$.  For all $i,j\in W$, just one of the sets $T_{ij}$ could be non-empty because $G^{26}$ is a forbidden subgraph. Moreover, $|T_{ij}|\leq1$ because $F_{1}$ and $E$ are forbidden subgraphs. Also, $\forall i \in W$, $|T_{i}|\leq2$ since $F_{6}, E, G^{16}$ and $S_{4}$ are forbidden subgraphs. Without loss of generality, let $T_{12}=\{x\}$. For all $i,j,k\in V''$, since $G^{23}, E, C_{4}, G^{24}$ are forbidden subgraphs,we get the following results; $T_{x}=T_{1}=T_{2}=T_{xij}=T_{xijk}=T_{x3}=T_{x4}=\emptyset$, $|T_{x1}|\leq1$ and $|T_{x2}|\leq1$. Hence the maximum value of the $t(H)$ could be equal to $\left(
               \begin{array}{c}
                 w(H) \\
                 3 \\
               \end{array}
             \right)+w(H)
$. When $w(H)<p$ and $w(H)\leq6$, we get $[\left(
               \begin{array}{c}
                 w(H) \\
                 3 \\
               \end{array}
             \right)+w(H)]<\left(
                                   \begin{array}{c}
                                     p \\
                                     3 \\
                                   \end{array}
                                 \right)
             $. Thus we have $T_{ij}=\emptyset$ for all $i,j \in V''$. Even though the sets $T_{i}$ are non-empty, the maximum value of the $t(H)$ could be equal to at most $\left(
                                         \begin{array}{c}
                                           w(H) \\
                                           3 \\
                                         \end{array}
                                       \right)+w(H)
             $ and this value is again less than $\left(
                                                    \begin{array}{c}
                                                      p \\
                                                      3 \\
                                                    \end{array}
                                                  \right)
             $. So this is a contradiction.

Hence, we obtain the fact in our second claim that H is not a generalized line graph of any graph.

 This means that, our assumption, that is $G$ is a disconnected graph, is actually false. Hence $G$ must be a connected graph.
\end{prf}

We give the our main result in the following theorem.

\begin{thm}
The $Kite_{p}^{q}$ graph is determined by its adjacency spectrum for all p and q.

\end{thm}

\begin{prf}
This result could be obtained directly from combining Lemma 3.3. and Corollary 3.2.
\end{prf}

\textbf{Acknowledgement}

The authors are indebted to  Willem H. Haemers for his valuable comments of the manuscript which significantly guides
our original arguments.



\end{document}